\documentclass[11pt]{amsart}

\usepackage{amsmath,amsthm}
\usepackage{amssymb}
\usepackage{euscript}
\numberwithin{equation}{section}

\newcommand{\sqsp}{\renewcommand{\baselinestretch}{1.1}\tiny\normalsize}

\newtheorem{thm}{Theorem}[section]
\newtheorem{prop}[thm]{Proposition}

\theoremstyle{definition}

\newcommand{\bZ}{\mathbf{Z}}

\DeclareMathOperator{\Id}{Id}

\begin{document}

\title{Loop structures on the homotopy type of $S^3$ revisited}
\author{Donald Yau}
\address{Department of Mathematics, University of Illinois at Urbana-Champaign, 1409 W.\ Green Street, Urbana, IL 61801}
\email{dyau@math.uiuc.edu}


\maketitle

\sqsp
\section{Introduction}

In an attempt to understand Lie groups from a homotopy theory point of view, Rector suggested studying Lie groups through their classifying spaces.  Using $S^3$ as a test case, he proved in his pioneering paper \cite{rector1} that there are uncountably many homotopically distinct deloopings of $S^3$.  These deloopings form the so-called \emph{genus} of the classifying space $BS^3$.  To be more precise, for a nilpotent finite type space $X$, the \emph{genus} of $X$ is defined to be the set of homotopy types of nilpotent finite type spaces $Y$ such that the $p$-completions of $X$ and $Y$ are homotopy equivalent for each prime $p$ and also their rationalizations are homotopy equivalent.  One often ignores the difference between a homotopy type and a space with that homotopy type when considering genus.

Rector actually provided a complete list of classification invariants, defined by using integral and mod $p$ cohomology, for the genus of $BS^3$, which we now call the Rector invariants.  Briefly, the Rector invariants of a space $X$ in the genus of $BS^3$ are signs, $(X/p) \in \lbrace \pm 1 \rbrace$, one for each prime $p$.  Two spaces in the genus of $BS^3$ are homotopy equivalent if and only if they have the same corresponding Rector invariants, and any such sequence does occur for some space.  Also, $BS^3$ itself has $1$ as its Rector invariant for each prime.

Generalizing the approach used by Rector, M\o ller \cite{moller} showed that this property of having a huge genus is not restricted to $BS^3$.   In fact, M\o ller proved that whenever $G$ is a compact connected non-abelian Lie group, the genus of its classifying space $BG$ is uncountably large.  There is, however, no explicit list of classification invariants, like the Rector invariants, in this general setting.  A remarkable result of Notbohm \cite{notbohm} showed that $K$-theory ring together with the $\lambda$-operations classify the genus of $BG$, provided $G$ is a simply-connected compact Lie group.  In other words, given two spaces in the genus of $BG$, they are homotopy equivalent if and only if their $K$-theory are isomorphic as $\lambda$-rings.  Notbohm's proof does not involve computing the $K$-theory of these spaces, and consequently we do not know how they are mutually non-isomorphic.  This prompts the question as to \emph{how} these uncountably many $\lambda$-rings, the $K$-theory of the spaces in the genus of a given classifying space $BG$, are mutually non-isomorphic.

Ideally, one would like to compute explicitly these $K$-theory $\lambda$-rings, at least partially.  Then one uses this knowledge to show that these $\lambda$-rings are mutually non-isomorphic by doing pure algebra.  Any work in this direction would shed new light into the results of Notbohm, M\o ller, and Rector mentioned above.

The purpose of this note is to do just that in the case $BS^3$.  In fact, the author has already shown in \cite{yau} that spaces in the genus of $BS^3$ all have isomorphic $K$-theory rings, namely, the power series ring $\bZ \lbrack \lbrack x \rbrack \rbrack$ with $x$ in filtration $4$.  (The result in \cite{yau} actually shows, more generally, that the underlying $K$-theory ring is constant throughout the genus of $BG$ as long as $BG$ has torsionfree integral cohomology.)  We will use both $KO$-theory and $K$-theory, since in order to define $(X/2)$ one has to consider either symplectic or orthogonal $K$-theory.  What we actually will do is to first
observe that the Rector invariants show up in the Adams operations in (both unitary and orthogonal) $K$-theory.  Then we show that these $\lambda$-rings are pairwise non-isomorphic by doing simple algebra about power series.

In addition to providing some new insights about the spaces in the genus of $BS^3$, the author also hopes that this note will attract more attention to the general problem of determining classification invariants for the genus of a classifying space $BG$.  Such classification invariants (if available) together with the author's result in \cite{yau} will probably allow one to show explicitly \emph{how} the $K$-theory $\lambda$-rings of the spaces in the genus of $BG$ differ from one another.

In the next section, we will first recall some basic definitions concerning $\lambda$-rings and Adams operations.  Those readers who are familiar with $\lambda$-rings can safely skip this section.  Our main results are presented in the last section, which begins with the construction of the Rector invariants.


\section{$\lambda$-rings and Adams operations}
\label{sec:lambda-ring}

In this section, we recall the definitions of a (filtered) $\lambda$-ring and of Adams operations.  The reader is referred to Atiyah and Tall \cite{at} or Knutson \cite{knutson} for more information about $\lambda$-rings.

\subsection{$\lambda$-rings}
A $\lambda$-\emph{ring} is a commutative ring $R$ with unit together with functions 
   \[
   \lambda^i \colon R ~\to~ R \quad (i = 0, 1, 2, \ldots)
   \]
such that for any elements $r$ and $s$ in $R$, one has
\begin{itemize}
\item $\lambda^0(r) = 1$
\item $\lambda^1(r) = r$, $\lambda^n(1) = 0$ for all $n > 1$
\item $\lambda^n(r + s) = \sum_{i = 0}^n\, \lambda^i(r)\lambda^{n-i}(s)$
\item $\lambda^n(rs) = P_n(\lambda^1(r), \ldots, \lambda^n(r);\lambda^1(s), \ldots, \lambda^n(s))$
\item $\lambda^n(\lambda^m(r)) = P_{n,m}(\lambda^1(r), \ldots, \lambda^{nm}(r))$.
\end{itemize}
Here the $P_n$ and $P_{n,m}$ are certain universal polynomials with integer coefficients.  Note that in the literature (for example, Atiyah and Tall \cite{at}) the terminology ``special'' $\lambda$-ring is used.

A $\lambda$-ring map $f \colon R \to S$ between two $\lambda$-rings is a ring map between the underlying rings which respects the operations $\lambda^i$: $f\lambda^i = \lambda^i f$ $(i \geq 0)$.

\subsection{Filtered $\lambda$-rings}
By a filtered ring we mean a commutative ring $R$ with unit together with a decreasing sequence of ideals, 
   \[
   R = I^0 \supset I^1 \supset I^2 \supset \cdots
   \]
such that $I^nI^m  \subset I^{n+m}$ for any $n,m$.  A map of filtered rings is a map of rings that respects the filtrations, in the sense that the image of the $n$th filtration is inside the $n$th filtration of the target.  (Note that this is a stronger condition than just saying that the ring map is continuous under the filtration topology.)

A \emph{filtered} $\lambda$-\emph{ring} is a filtered ring $R = (R, \lbrace I^n \rbrace)$ which is also a $\lambda$-ring such that the ideals $I^n$ are all closed under the operations $\lambda^i$ $(i > 0)$.

A filtered $\lambda$-ring map is a $\lambda$-ring map which is also a filtered ring map.

\subsection{Adams operations}
Given a $\lambda$-ring $R$, the \emph{Adams operations} 
   \[
   \psi^k \colon R ~\to~ R
   \]
for $k = 1, 2, \ldots$ are defined inductively by the Newton formula:
   \[
   \psi^k(a) - \lambda^1(a) \psi^{k-1}(a) + \cdots + (-1)^{k-1}\lambda^{k-1}(a) \psi^1(a) = (-1)^{k-1}k \lambda^k(a).
   \]
The Adams operations satisfy the following properties.
   \begin{enumerate}
   \item All the $\psi^k \colon R \to R$ are $\lambda$-ring maps.
   \item $\psi^1 = \Id$ and $\psi^k \psi^l = \psi^{kl}$ for any $k, l \geq 1$.
   \item $\psi^p(a) \equiv a^p$ (mod $pR$) for each prime $p$ and element $a$ in $R$.
   \end{enumerate}

If $R$ is a filtered $\lambda$-ring, then the Adams operations are all filtered $\lambda$-ring maps.  Observe from the second condition that the $\psi^p$ for $p$ primes determine all the Adams operations.


\section{Rector invariants and the genus of $BS^3$}

Let $X$ be a space in the genus of $BS^3$.  Using the proof of the Main Theorem in \cite{yau} with $KO$-theory instead of $K$-theory, one can show that
   \[
   KO^*(X) \cong KO^* \lbrack \lbrack x \rbrack \rbrack
   \]
as filtered rings, where $x \in KO^4_4(X)$ is a representative of an integral generator $x_4 \in H^4(X;\bZ) = E_2^{4,\, 0}$ in the $KO^*$-Atiyah-Hirzebruch spectral sequence for $X$.  Here $KO^a_b(X)$ denotes the subgroup of $KO^a(X)$ consisting of elements $u$ which restrict to $0$ under the natural map 
   \[
   KO^a(X) ~\to~ KO^a(X_{b-1}),  
   \]
where $X_{b-1}$ denotes the $(b-1)$-skeleton of $X$.  Such an element $u$ is said to be in \emph{degree} $a$ and \emph{filtration} $b$.

Now we recall the relevant notations, definitions, and results regarding Rector's classification of the genus of $BS^3$ \cite{rector1}.  Let $\xi \in \pi_{-4} KO$ and $b_R \in \pi_{-8} KO$ be the generators so that $\xi^2 = 4b_R$.  Rector observed in \cite[\S 4]{rector2} that there exists an integer $a$, depending on the choice of the representative $x$, such that the following three statements hold.
   \begin{enumerate}
   \item $\psi^2(\xi x) = 4 \xi x + 2a b_R x^2$ (mod $KO^0_9(X)$).
   \item The integer $a$ is well-defined (mod $24$).  This means that if $x^\prime \in KO^4_4(X)$ is another representative of $x_4$ with corresponding integer $a^\prime$, then $a \equiv a^\prime$ (mod $24$).  Furthermore, if $x_4$ is replaced with $-x_4$, then $a$ will be replaced with $-a$.  We can (and we will) therefore write $a(X)$ for $a$.
   \item $a(X) \equiv \pm 1, \pm 5, \pm 7,$ or $\pm 11$ (mod $24$).  
   \end{enumerate}
The last condition above follows from the examples constructed by Rector in \cite[\S 5]{rector2} and James' result \cite{james} which says that there are precisely eight homotopy classes of homotopy-associative multiplications on $S^3$.  These eight classes can be divided into four pairs with each pair consisting of a homotopy class of multiplication and its inverse.

Rector's invariant $(X/p)$ for $p$ an odd prime is defined as follows \cite{rector1}.  In the mod $p$ Steenrod algebra, the Adem relation $P^1 P^1 = 2P^2$ implies that 
   \[
   P^1 \, \overline{x}_4 ~=~ \pm \, 2\overline{x}_4^{(p+1)/2}
   \]
in $H^*(X;\bZ/p)$, where $\overline{x}_4$ is the mod $p$ reduction of the integral generator $x_4$.  Then $(X/p) \in \{ \pm 1\}$ is defined as the sign on the right-hand side of this equation.  Since replacing $x_4$ by $-x_4$, the other integral generator, will change all the signs $(X/p)$ simultaneously, one also has to specify an orientation of the integral generator in order that the $(X/p)$ be well-defined.

The Rector invariant $(X/2)$ and a canonical choice of an orientation of the integral generator are given as follows.  Using the (mod $24$) integer $a(X)$, define
   \begin{equation}
   \label{eq0:BS3}
   ((X/2), (X/3)) \buildrel \text{def} \over = 
   \begin{cases}
   (1,1) & \text{if }a(X) \equiv \pm 1 \text{ mod } 24;\\
   (1,-1) & \text{if }a(X) \equiv \pm 5 \text{ mod } 24;\\
   (-1,1) & \text{if }a(X) \equiv \pm 7 \text{ mod } 24;\\ 
   (-1,-1) & \text{if }a(X) \equiv \pm 11 \text{ mod } 24.
   \end{cases}
   \end{equation}
The orientation of the integral generator is then chosen so that $(X/3)$ is as given in \eqref{eq0:BS3}.  This definition of Rector's invariants coincides with the original one (cf.\ \cite[\S 9]{mcgibbon}).  In particular, specifying the Rector invariants of a space $X$ is equivalent to specifying the (mod $24$) integer $a(X)$ and the $(X/p)$ for $p$ primes, $p > 2$.  Observe from this discussion that $a(X)$ appears in the Adams operation $\psi^2$ in $KO$-theory.

Now we can recall Rector's classification theorem of the genus of $BS^3$ \cite{rector1}.

\medskip
\begin{thm}[Rector]
\label{thm:rector}
The $(X/p)$ for $p$ primes provide a complete list of classification invariants for the genus of $BS^3$.  Any combination of values of the $(X/p)$ can occur.  If $X$ is $BS^3$, then $(X/p) = 1$ for every prime $p$.
\end{thm}
\medskip

It should be remarked that Rector also proved that $BS^3$ can be distinguished from the rest of the spaces in its genus by considering certain maps, the so-called ``maximal torus," from infinite complex projective space.   However, we will not need to use that result here.

We have seen above how the (mod 24) integer $a(X)$ appears in the Adams operation $\psi^2$ in $KO$-theory.  Now we will see how the Rector invariants $(X/p)$ for $p$ odd appear in the Adams operations $\psi^p$ in $K$-theory.

\medskip
\begin{prop}
\label{prop:X/p}
Let $p$ be a fixed odd prime and let $X$ be a space in the genus of $BS^3$, so that $K^*(X)$ can be written as $K^*\lbrack \lbrack u \rbrack \rbrack$ for some representative $u \in K^4_4(X)$ of the integral generator $x \in H^4(X;\bZ) = E_2^{4,0}$ in the $K^*$-Atiyah-Hirzebruch spectral sequence.  Then there exist some elements $w \in K^0_{2p+3}(X)$ and $x_0 \in K^0_4(X)$ such that
   \[
   \psi^p(b^2 u) ~=~ (b^2 u)^p + 2(X/p)\, p\, (b^2 u)^{(p+1)/2} + p\, w + p^2 x_0,
   \]
where $b \in \pi_{-2}K$ is the Bott element.
\end{prop}
\medskip

\begin{proof}
To see this, note that since $b^2 u \in K^0_4(X)$ is in filtration $4$, it follows from Atiyah's theorem \cite[Proposition 5.6]{atiyah} that
   \[
   \psi^p(b^2 u) ~=~ (b^2 u)^p + p\, x_1 + p^ 2x_0
   \]
for some $x_i \in K^0_{4+2i(p-1)}(X)$ $(i = 0, 1)$.  Moreover, one has 
   \[
   \overline{x}_1 ~=~ P^1\, \overline{b^2 u}, 
   \]
where $\overline{z}$ is the mod $p$ reduction of an element $z$ and $P^1$ is the Steenrod operation of degree $2(p-1)$ in mod $p$ cohomology.  Thus, to prove the Proposition it is enough to show that
   \begin{equation}
   \label{eq4:BS3}
   x_1 ~=~ 2(X/p)(b^2 u)^{(p+1)/2} + w + p\, z
   \end{equation}
for some $w \in K^0_{2p+3}(X)$ and some $z \in K^0_{2p+2}(X)$.   Now in $H^*(X;\bZ) \otimes \bZ/p = H^*(X;\bZ/p)$ we have
   \[
   \begin{split}
   \overline{x}_1 & ~=~ P^1\, \overline{b^2 u} \\
               & ~=~ P^1\, \overline{x}     \\
               & ~=~ 2(X/p)\, \overline{x}^{(p+1)/2} \\
               & ~=~ 2(X/p)\, \overline{b^2 u}^{(p+1)/2}.
   \end{split}
   \]
The Proposition follows from this because the mod $p$ cohomology of $X$ is the mod $p$ associated graded ring of its $K$-theory.
\end{proof}

Equipped with the knowledge of how the Rector invariants appear in $K$-theory, we now show that the $K$-theory of the spaces in the genus of $BS^3$ can be distinguished by doing some rather elementary algebra.

\medskip
\begin{thm}
Let $X$ and $Y$ be spaces in the genus of $BS^3$.  If there exists a filtered $\lambda$-ring isomorphism $\sigma \colon KO^*(X) \xrightarrow{\cong} KO^*(Y)$, then 
   \[
   a(X) ~\equiv~ \pm \, a(Y) \pmod{24}
   \]
and 
   \[
   (X/p) ~=~ (Y/p)
   \]
for each odd prime $p$.  In this case, by Rector's classification theorem, it follows that $X$ and $Y$ have the same homotopy type.
\end{thm}
\medskip

\begin{proof}
As explained above, $KO^*(X) = KO^*\lbrack \lbrack x \rbrack \rbrack$ and $KO^*(Y) = KO^*\lbrack \lbrack y \rbrack \rbrack$ with $x \in KO^4_4(X)$ and $y \in KO^4_4(Y)$ representing, respectively, the integral generators $x_4 \in H^4(X;\bZ)$ and $y_4 \in H^4(Y;\bZ)$.

Since $\sigma$ is a ring isomorphism, we have 
   \[
   \sigma(\xi x) ~=~ \varepsilon \xi y + \sigma_2 b_R y^2 \quad (\text{mod } KO^0_9(Y)) 
  \]
for some $\varepsilon \in \{\pm 1\}$ and integer $\sigma_2$.  Computing modulo $KO^0_9(X)$ we have
   \[
   \begin{split}
    4\sigma(b_Rx^2) & ~=~ \sigma(\xi x)^2 \\
                    & ~=~ \xi^2 y^2       \\
                    & ~=~ 4b_R y^2.
   \end{split}
   \]
Therefore, one has 
   \[
   \sigma(b_Rx^2) ~=~ b_Ry^2 \quad (\text{mod } KO^0_9(Y)).  
   \]
First we claim that there is an equality
   \begin{equation}
   \label{eq1:BS3}
   a(X) ~=~ 6 \sigma_2 + \varepsilon a(Y).
   \end{equation}
To prove \eqref{eq1:BS3} we will compute both sides of the equality
   \[
   \sigma \psi^2(\xi x) ~=~ \psi^2 \sigma (\xi x) \quad  (\text{mod } KO^0_9(Y)).  
   \]
Working modulo $KO^0_9(Y)$ we have, on the one hand,
   \begin{eqnarray*}
   \sigma \psi^2(\xi x) 
   & ~=~ & \sigma (4 \xi x + 2a(X)b_R x^2) \\
   & ~=~ & 4(\varepsilon \xi y + \sigma_2 b_R y^2) + 2a(X)b_Ry^2 \\
   & ~=~ & 4\varepsilon \xi y + (4 \sigma_2 + 2a(X))b_Ry^2.  
   \end{eqnarray*}
On the other hand, still working modulo $KO^0_9(Y)$, we have
   \begin{eqnarray*}
   \psi^2 \sigma (\xi x) 
   & ~=~ & \varepsilon\psi^2(\xi y) + \sigma_2 \psi^2(b_Ry^2) \\
   & ~=~ & \varepsilon(4\xi y + 2a(Y)b_Ry^2) + \sigma_2(2^4 b_Ry^2) \\
   & ~=~ & 4 \varepsilon \xi y + (16\sigma_2 + 2\varepsilon a(Y))b_Ry^2.
   \end{eqnarray*}   
Equation \eqref{eq1:BS3} now follows by equating the coefficients of $b_Ry^2$.

In view of \eqref{eq1:BS3}, to finish the proof of the assertion about the mod $24$ integer $a$, it is enough to establish
   \begin{equation}
   \label{eq2:BS3}
   \sigma_2 ~\equiv~ 0 \quad (\text{mod }4).
   \end{equation}
To prove \eqref{eq2:BS3}, note that since $\sigma$ is a $KO^*$-module map, we have 
   \[
   \xi \sigma(x) ~=~ \sigma(\xi x).  
   \]
Since $\sigma$ is a ring isomorphism, we also have 
   \[
   \sigma(x) ~=~ \varepsilon^\prime y + \sigma_2^\prime \xi y^2 \quad (\text{mod } KO^4_9(Y)) 
   \]
for some $\varepsilon^\prime \in \{\pm 1\}$ and integer $\sigma_2^\prime$.  Therefore, working modulo $KO^0_9(Y)$ we have
   \begin{eqnarray*}
   \xi \sigma(x) 
   & ~=~ & \varepsilon^\prime \xi y + \sigma_2^\prime \xi^2 y^2 \\
   & ~=~ & \varepsilon^\prime \xi y + 4\sigma_2^\prime b_R y^2 \\
   & ~=~ & \varepsilon \xi y + \sigma_2 b_R y^2.
   \end{eqnarray*}
In particular, by equating the coefficients of $b_R y^2$ we obtain 
   \[
   \sigma_2 ~=~ 4 \sigma_2^\prime, 
   \]
thereby proving \eqref{eq2:BS3}.

This finishes the proof of the assertion that $a(X) \equiv \pm \, a(Y) \pmod{24}$.  We still need to prove the assertion about the Rector invariants for odd primes $p$.

Let, then, $p$ be a fixed odd prime.  Just as before, we have that
   \[
   K^*(X) ~\cong~ K^*\lbrack \lbrack u_x \rbrack \rbrack
   \] 
with $u_x \in K^4_4(X)$ a representative of the integral generator $x_4 \in H^4(X;\bZ) = E_2^{4,0}$ in the $K^*$-Atiyah-Hirzebruch spectral sequence.  Moreover, we may choose $u_x$ so that 
   \[
   c(x) ~=~ u_x, 
   \]
where 
   \[
   c \colon KO^*(X) ~\to~ K^*(X)
   \] 
is the complexification map.  Similar remarks apply to $Y$ so that 
   \[
   K^*(Y) ~\cong~ K^*\lbrack \lbrack u_y \rbrack \rbrack.
   \]

The $\lambda$-ring isomorphism $\sigma$ induces via $c$ a $\lambda$-ring isomorphism  
   \[
   \sigma_c \colon K^*(X) ~\xrightarrow{\cong}~ K^*(Y).  
   \]
By composing $\sigma_c$ with a suitable $\lambda$-ring automorphism of $K^*(Y)$ if necessary, we obtain a $\lambda$-ring isomorphism 
   \[
   \alpha \colon K^*(X) ~\xrightarrow{\cong}~ K^*(Y) 
   \]
with the property that
   \begin{equation}
   \label{eq5:BS3}
   \alpha(b^2 u_x) ~=~ b^2 u_y + \text{ higher terms in }b^2 u_y.
   \end{equation}
Using Proposition \ref{prop:X/p} and \eqref{eq5:BS3}, one infers that
   \begin{equation}
   \label{eq6:BS3}
   \alpha \psi^p(b^2 u_x) ~=~ 2(X/p) \, p \, (b^2 u_y)^{(p+1)/2} \quad    (\text{mod } K^0_{2p+3}(Y) \text{ and } p^2)
   \end{equation}
and
   \begin{equation}
   \label{eq7:BS3}
\psi^p \alpha(b^2 u_x) ~=~ 2(Y/p) \, p \, (b^2 u_y)^{(p+1)/2} \quad (\text{mod } K^0_{2p+3}(Y) \text{ and } p^2).
   \end{equation}
Since $\alpha \psi^p = \psi^p \alpha$ it follows from \eqref{eq6:BS3} and \eqref{eq7:BS3} that
   \[
   2(X/p) \, p ~\equiv~ 2(Y/p) \, p \quad (\text{mod } p^2),
   \]
or, equivalently,
   \[
   2(X/p) ~\equiv~ 2(Y/p) \quad (\text{mod } p).
   \]
But $p$ is odd, and so 
   \[
   (X/p) ~\equiv~ (Y/p) \quad (\text{mod } p).  
   \]
Hence $(X/p) = (Y/p)$, as desired.

This finishes the proof of the Theorem.
\end{proof}

\section{Acknowledgment}
The results in this paper constitute part of the author's 2002 Ph.D.\ thesis at MIT.  The author expresses his sincerest gratitude to his advisor Professor Haynes Miller for his guidance, mathematical and otherwise.


\end{document}